\documentclass[10 pt]{amsart}
\usepackage[applemac]{inputenc}

\usepackage[dvipsnames,svgnames,x11names,hyperref]{xcolor}
\usepackage[hyphens]{url}
\usepackage[colorlinks=true,linkcolor=Maroon,citecolor=blue,urlcolor=blue,hypertexnames=false,linktocpage]{hyperref}
\usepackage{bookmark}
\usepackage{amsthm,thmtools,amssymb,amsmath,amscd}

\usepackage{fancyhdr}
\usepackage{esint}
\usepackage{enumerate}

\usepackage{pictexwd,dcpic}
\usepackage{graphicx}
\usepackage{caption}
\swapnumbers
\declaretheorem[name=Theorem,numberwithin=section]{thm}
\declaretheorem[name=Remark,style=remark,sibling=thm]{rem}
\declaretheorem[name=Lemma,sibling=thm]{lemma}

\declaretheorem[name=Definition,style=definition,sibling=thm]{defn}
\declaretheorem[name=Corollary,sibling=thm]{cor}
\declaretheorem[name=Assumption,style=definition,sibling=thm]{assum}

\numberwithin{equation}{section}

\usepackage{cleveref}
\crefname{lemma}{Lemma}{Lemmata}
\crefname{prop}{Proposition}{Propositions}
\crefname{thm}{Theorem}{Theorems}
\crefname{cor}{Corollary}{Corollaries}
\crefname{defn}{Definition}{Definitions}
\crefname{example}{Example}{Examples}
\crefname{rem}{Remark}{Remarks}
\crefname{assum}{Assumption}{Assumptions}
\crefname{nota}{Notation}{Notation}
\usepackage{autonum}

\newcommand{\wt}{\widetilde}
\newcommand{\wh}{\widehat}

\newcommand{\cn}{\colon}
\newcommand{\sub}{\subset}
\newcommand{\ov}{\overline}

\newcommand{\bbR}{\mathbb{R}}

\newcommand{\bbS}{\mathbb{S}}
\newcommand{\bbH}{\mathbb{H}}

\newcommand{\bbE}{\mathbb{E}}

\newcommand{\8}{\infty}

\newcommand{\al}{\alpha}
\newcommand{\be}{\beta}

\newcommand{\de}{\delta}
\newcommand{\ep}{\epsilon}
\newcommand{\ka}{\kappa}

\newcommand{\si}{\sigma}

\newcommand{\ph}{\phi}

\newcommand{\vt}{\vartheta}

\newcommand{\De}{\Delta}
\newcommand{\Ga}{\Gamma}
\newcommand{\Th}{\Theta}

\newcommand{\cF}{\mathcal{F}}

\newcommand{\cS}{\mathcal{S}}

\newcommand{\cP}{\mathcal{P}}

\newcommand{\del}{\partial}
\newcommand{\n}{\nabla}

\newcommand{\fa}{\forall}
\newcommand{\rt}{\sqrt}

\newcommand{\fr}[2]{\frac{#1}{#2}}
\newcommand{\tfr}[2]{\tfrac{#1}{#2}}
\newcommand{\x}{\times}

\DeclareMathOperator{\osc}{osc}
\DeclareMathOperator{\const}{const}

\DeclareMathOperator{\Rm}{Rm}
\DeclareMathOperator{\Rc}{Rc}

\newcommand{\pf}[1]{\begin{proof}#1 \end{proof}}
\newcommand{\eq}[1]{\begin{equation}\begin{alignedat}{2} #1 \end{alignedat}\end{equation}}

\newcommand{\br}[1]{\left(#1\right)}
\newcommand{\abs}[1]{\lvert #1\rvert}
\newcommand{\enum}[1]{\begin{enumerate}[(i)] #1 \end{enumerate}}
\newcommand{\enu}[1]{\begin{enumerate}[(a)] #1 \end{enumerate}}

\newcommand{\ra}{\rightarrow}
\newcommand{\hra}{\hookrightarrow}

\newcommand{\mt}{\mapsto}

\newcommand{\hp}{\hphantom}
\newcommand{\q}{\quad}

\begin{document}
\title[Locally constrained flows in warped spaces]{Minkowski inequalities and constrained inverse curvature flows in warped spaces}

\date{\today}
\thanks{Funded by the "Deutsche Forschungsgemeinschaft" (DFG, German research foundation); Project "Quermassintegral preserving local curvature flows"; Grant number SCHE 1879/3-1.}
\author[J. Scheuer]{Julian Scheuer}
\address{Department of Mathematics, Columbia University
New York, NY 10027, USA}
\email{\href{mailto:jss2291@columbia.edu}{jss2291@columbia.edu}; \href{mailto:julian.scheuer@math.uni-freiburg.de}{julian.scheuer@math.uni-freiburg.de}}

\begin{abstract}
This paper deals with locally constrained inverse curvature flows in a broad class of Riemannian warped spaces. For a certain class of such flows we prove long time existence and smooth convergence to a radial coordinate slice. In the case of two-dimensional surfaces and a suitable speed, these flows enjoy two monotone quantities. In such cases new Minkowski type inequalities are the consequence.
In higher dimensions we use the inverse mean curvature flow to obtain new Minkowski inequalities when the ambient radial Ricci curvature is constantly negative.
\end{abstract}

\maketitle

\section{Introduction}		

The objectives of this paper are threefold. First we want to continue the investigation of the so-called {\it{locally constrained inverse curvature flows}}. These are hypersurface variations of the form
\eq{\label{intro:flow}\del_{t}{x}=-\cF(u,s,\ka)\nu,}
where $x$
is a smooth family of embeddings of a smooth compact manifold into an ambient Riemannian manifold 
\eq{\label{warped}N=(a,b)\x \cS_{0},\q \bar g=dr^{2}+\vt^{2}(r)\si.} 
$(\cS_{0},\si)$ is a compact Riemannian manifold of dimension $n\geq 2$ and $\vt\in C^{\8}([a,b))$.
The slices $M_{t}=x(t,\cS_{0})$ are graphical over $\cS_{0}$ with graph function $u$, support function $s$, principal curvatures $\ka$ and outward unit normal $\nu$.  In this paper we investigate flows of the form
\eq{\del_{t}{x}=\br{\fr{\vt'(u)}{F(\ka)}-s}\nu}
and prove convergence to a radial slice $\{r=\const\}$ under various assumptions on $F$ and $N$.

The second objective is to apply this result in case $n=2$ with the particular choice
\eq{F=\fr{H_{2}}{H_{1}}}
 in order to prove new geometric inequalities. Here $H_{k}$ is the $k$-th normalized elementary symmetric polynomial of the principal curvatures.
 
 Finally, we accompany these results by some new Minkowski inequalities  in higher dimension and for some ambient spaces of non-constant curvature. These are consequences of the inverse mean curvature flow.

We state the main results, after imposing some general assumptions and fixing some notation.

\begin{assum}\label{N}
For $n\geq 2$ let $(\cS_{0},\si)$ be a compact, $n$-dimensional Riemannian manifold, $a<b$ real numbers and $\vt\in C^{\8}([a,b))$.
We assume that the warped product space 
\eq{N=(a,b)\x\cS_{0},\q \bar g=dr^{2}+\vt^{2}(r)\si}
 satisfies the following assumptions:
\enum{
\item \label{N-1}$\vt'>0$,
\item \label{N-2}Either of the following conditions hold:
	\enu{ 
	\item $\vt''\geq 0$ 
	\item  $\vt''\leq 0$ and 
	\eq{\label{N-3}\del_{r}\br{\fr{\vt''}{\vt}}\leq 0.}
	}  
Furthermore we denote by $\wh{\Rc}$ the Ricci curvature of the metric $\si$ and geometric quantities of the ambient space $N$ are furnished with an overbar, e.g. $\ov{\Rm}$, $\ov{\Rc}$ and $\bar\n$ for the Riemann tensor, the Ricci tensor and the Levi-Civita connection respectively.
}
\end{assum}

\begin{rem}
Note that \eqref{N-3} says that the Ricci curvature of the ambient space is non-decreasing in radial directions.
\end{rem}

\begin{defn}\label{M}
Let $N$ be as in \Cref{N} and $M\sub N$ be a graph over $\cS_{0}$, i.e.
\eq{M=\{(u(y),y)\cn y\in \cS_{0}\}.} We define the {\it{region enclosed by $M$}} to be the set
\eq{\hat M=\{(r,y)\cn a\leq r\leq u(y),~y\in \cS_{0}\}.}
We denote by $\abs{M}$ the surface area of $M$ and by $\abs{\hat M}$ the volume of the region $\hat M$.
We call $M$ {\it{strictly convex}}, if all principal curvatures have a strict sign everywhere.
\end{defn}

\begin{assum}\label{F}
Let 
\eq{\Ga_{+}=\{\ka\in \bbR^{n}\cn \ka_{i}>0\q\fa 1\leq i\leq n\}.}
We suppose that $F\in C^{\8}(\Ga_{+})$ is symmetric, positive, strictly monotone, homogeneous of degree $1$, concave and satisfies
\eq{F_{|\del\Ga_{+}}=0,\q F(1,\dots,1)=1.} 
\end{assum}

Here is the main result concerning the curvature flow.

\begin{thm}\label{main:flow}
Let $N$ and $F$ satisfy \Cref{N,,F} respectively. Let
\eq{x_{0}\cn M_{0}\hra N}
 be the embedding of a strictly convex hypersurface, given as a graph over $\cS_{0}$,
\eq{M_{0}=\{(u_{0}(y),y)\cn y\in \cS_{0}\}.}
Then there exists a unique immortal solution
\eq{x\cn [0,\8)\x M_{0}\ra N}
that satisfies the parabolic Cauchy-problem
\eq{\label{flow}\del_{t}x&=\br{\fr{\vt'(u)}{F}-s}\nu\\
		x(0,\cdot)&=x_{0}.}
The embeddings $x(t,\cdot)$ converge smoothly to an embedding of a slice $\{r=\const\}$.
\end{thm}

Flows of the kind \eqref{flow} originated from an idea due to Guan and Li \cite{GuanLi:/2015}. At first they investigated the mean curvature type flow
\eq{\label{GL}\del_{t}{x}=(\vt'-sH_{1})\nu}
in the simply connected spaceforms of constant curvature and noticed that \eqref{GL} preserves the enclosed volume and decreases the surface area. This way they reproved the isoperimetric inequality in spaceforms for starshaped hypersurfaces. This flow was later transferred to more general warped product spaces by the same authors and Wang \cite{GuanLiWang:/2019}. The Lorentzian version is treated in \cite{LambertScheuer:10/2019}.

 It is easy to show that the general Minkowski identities
\eq{\int_{M}\vt'H_{k}=\int_{M}sH_{k+1}}
imply that the flow \eqref{flow} with
\eq{\label{Hk}F=\fr{H_{k}}{H_{k-1}}}
preserves certain quermassintegrals for starshaped hypersurfaces in spaceforms and decreases/increases the others. See the survey \cite{GuanLi:/2019} for more details. Hence it has been tempting to prove that the flow \eqref{flow} converges to a totally umbilic hypersurface in order to get new results in integral geometry. Unfortunately this flow has refused to release the results it originally promised, at least so far. There are some partial results in spaceforms. Denote by $\bbE^{n+1}$ the Euclidean space, by $\bbS_{+}^{n+1}$ the hemisphere, by $\bbH^{n+1}$ the hyperbolic space and by $\bbS^{n,1}_{+}$ the upper branch of the Lorentzian de Sitter space of respective dimension $n+1$. To summarize the known results, let $F$ be given by \eqref{Hk}. Then, \eqref{flow} starting from a starshaped hypersurface with $F>0$ converges to a totally umbilic hypersurface provided
\enum{
\item $N=\bbE^{n+1}$. This is a trivial case which already follows from the classical inverse curvature flows due to Gerhardt \cite{Gerhardt:/1990} and Urbas \cite{Urbas:/1990},
\item $N=\bbS^{n+1}$, $N=\bbH^{n+1}$, $k=n$, \cite{BrendleGuanLi:/,GuanLi:/2019},
\item $N=\bbH^{n+1}$, $1\leq k\leq n$ and the initial hypersurface is $h$-convex, i.e. $\ka_{i}>1$, \cite{HuLiWei:02/2020},
\item $N=\bbH^{n+1}$, $1\leq k\leq n$ and the initial hypersurface satisfies an a priori gradient bound
\eq{\max_{x\in\bbS^{n}}\abs{\n\log\vt'(0,x)}^{2}\leq 12+3\min_{x\in\bbS^{n}}\vt(0,x)^{2},}
\cite{BrendleGuanLi:/,GuanLi:/2019},
\item $N=\bbS^{n,1}_{+}$, $k=1$, see \cite{Scheuer:09/2019}.
}
In each of these cases, one obtains a corresponding quermassintegral inequality. 
\Cref{main:flow} provides the first convergence result for a flow of type \eqref{flow} outside the constant curvature spaces for a class of $F$ that contains the case $F=H_{n}/H_{n-1}$.

As an application we obtain two new Minkowski type inequalities for surfaces in a certain class of warped spaces. 

Note that although we prove the following two results for $n=2$, we keep notation general, in order to show that the monotonicity properties of the flow with $F=H_{2}/H_{1}$ are also valid in higher dimension. The restriction to $n=2$ stems from the technical hurdle that $F=H_{2}/H_{1}$ only vanishes on the boundary of $\Ga_{+}$ if $n=2$ and hence only in this case we have a good convergence result.

\begin{thm}\label{thm:mink}
Suppose $n=2$ and in addition to \Cref{N} suppose that $N$ satisfies
\eq{\label{thm:mink-2}\wh{\Rc}\geq (n-1)(\vt'^{2}(r)-\vt''\vt(r))\si\q \fa r.}
Let $M\sub N$ be a strictly convex graph over $\cS_{0}$.
Then there holds
\eq{\label{thm:mink-1}\int_{M}H_{1}+\fr 1n\int_{\hat M}\ov\Rc(\del_{r},\del_{r})\geq \ph(\abs{M}),}
 where $\ph$ is the function that gives equality on the radial slices. If equality holds, then $M$ is totally umbilic. If the associated quadratic forms in \eqref{thm:mink-2} satisfy the strict inequality on nonzero vectors, then equality in \eqref{thm:mink-1} holds precisely on radial slices.
\end{thm}

\begin{rem}\label{rem:mink}
\enum{
\item Note that in case $\vt''\leq 0$, there holds
\eq{\del_{r}(\vt'^{2}-\vt''\vt)=-\vt^{2}\del_{r}\br{\fr{\vt''}{\vt}}\geq 0}
by assumption \eqref{N-3}. 
\item Note that \Cref{thm:mink} even holds in ambient spaces where the validity of the isoperimetric inequality is unclear, compare \cite[Sec.~6]{GuanLiWang:/2019}. In particular we do NOT assume that
\eq{\vt'^{2}-\vt''\vt\geq 0.} 
}
\end{rem}

Due to \Cref{rem:mink} it is also of interest to obtain a lower bound in terms of the volume of $\hat M$, as this can not be covered by an isoperimetric inequality. It is possible to deduce such an inequality under presence of a {\it{Heintze-Karcher type inequality.}}

\begin{thm}\label{thm:HK}
Suppose $n=2$ and in addition to \Cref{N} suppose that $\vt''\geq 0$ and
\eq{\label{thm:HK-3}\wh{\Rc}\geq (n-1)(\vt'^{2}(r)-\vt''\vt(r))\si\q \fa r.}
 Let $M\sub N$ be a strictly convex graph over $\cS_{0}$ and suppose for every such $M$ there holds
\eq{\label{thm:HK-2}\int_{M}\fr{\vt'}{H_{1}}\geq \int_{M}s }
and that equality implies total umbilicity.

Then
\eq{\label{thm:HK-1}\int_{M}H_{1}+\fr 1n\int_{\hat M}\ov\Rc(\del_{r},\del_{r})\geq \psi(\abs{\hat M}),}
where $\psi$ is the function which gives equality on the radial slices. If equality holds, then $M$ is totally umbilic. If the associated quadratic forms in \eqref{thm:HK-3} satisfy the strict inequality on nonzero vectors, then equality in \eqref{thm:HK-1} holds precisely on radial slices.
\end{thm}

\begin{rem}
In \cite[Thm.~3.11]{Brendle:06/2013}, inequality \eqref{thm:HK-2} was proved under the assumptions that $\vt''>0$ and that the ambient space is {\it{substatic}}, i.e.
\eq{\bar\De\vt'\bar g-\bar\n^{2}\vt'+\vt'\ov{\Rc}\geq 0.}
It is unclear however, whether this condition is necessary for the Heintze-Karcher inequality to hold.
\end{rem}

When it comes to higher dimensions, the flow \eqref{flow} with $F=H_{2}/H_{1}$ is not understood yet. However, if we restrict the ambient space further, it is possible to obtain a Minkowski inequality, which holds for $n\geq 2$, provided we impose a special structure of the ambient space in radial direction. It's proof does not rely on a locally constrained curvature flow, but on the standard inverse mean curvature flow. The idea on how to exploit its monotonicity properties comes from  \cite{BrendleGuanLi:/}, also see \cite{GuanLi:/2019}, in which the hyperbolic case is treated. The convergence of the inverse curvature flows in general warped products was proven in \cite{Scheuer:02/2019}, also see \cite{Zhou:04/2018}. However note that the latter work does not cover the required asymptotics in the ambient spaces we are considering here.

\begin{thm}\label{gen-mink}
In addition to \Cref{N} suppose that $(\cS_{0},\si)$ has non-negative sectional curvature. Suppose 
\eq{\vt(r)=\al\sinh(r)+\be\cosh(r),}
where $\al\geq\be\geq 0$ and one of those inequalities has to be strict and
\eq{\label{gen-mink-2}\wh\Rc\geq (n-1)(\al^{2}-\be^{2}) \si.}
Let $M\sub N$ be a strictly mean-convex graph over $\cS_{0}$. Then there holds
\eq{\label{gen-mink-1}\int_{M}H_{1}-\abs{\hat M}\geq \phi(\abs{M}),}
where $\phi$ is the function that gives equality on the radial slices.
If equality holds, then $M$ is totally umbilic. If the associated quadratic forms in \eqref{gen-mink-2} satisfy the strict inequality on nonzero vectors, then equality in \eqref{gen-mink-1} holds precisely on radial slices.
\end{thm}

Let us put the results in \Cref{thm:mink,,thm:HK,,gen-mink} into some historical context.
The functional
\eq{W_{2}(M):=\int_{M}H_{1}+\fr{1}{n}\int_{\hat M}\ov{\Rc}(\del_{r},\del_{r})}
plays a significant role in hypersurface theory and especially in the theory of convex bodies. In the Euclidean space and up to a dimensional constant, it arises from the Taylor expansion of volume with respect to outward geodesic variations of a compact domain $\hat M$ with smooth boundary $M$, which is compressed into the beautiful {\it{Steiner formula}} for convex bodies in Euclidean space \cite{Steiner:/2013},
\eq{\abs{\hat M_{\ep}}=\sum_{k=0}^{n+1}c_{n,k}W_{k}(M)\ep^{k},}
where $\hat M_{\ep}$ is the $\ep$-parallel body of $\hat M$. There are many related formulae for domains of the hyperbolic and spherical spaces, \cite{Allendoerfer:/1948,Santalo:06/1950}, and also see \cite{GallegoSolanes:/2005,GaoHugSchneider:/2001} for good introductions. The additional Ricci term has reasons stemming from a particular geometric interpretation of the Steiner coefficients, which requires the additional Ricci term when transferred to other ambient spaces. Regardless of the ambient space however, the Minkowski inequality provides a {\it{convexity estimate}} for the function
\eq{\ep\mt \abs{\hat M_{\ep}}}
and estimates its second derivative at $\ep=0$ from below by its value and its first derivative at $\ep=0$. As such, it makes a statement of volume growth and hence is of interest in Riemannian geometry. There has been immense effort in the past to obtain Minkowski inequalities, even for non-convex hypersurfaces. In the Euclidean space this was accomplished for convex bodies in \cite{Minkowski:/1903} and for starshaped and mean-convex hypersurfaces in \cite{GuanLi:08/2009}. It is open until today, whether the starshapedness can be dropped here. It can be shown however, that the result is also true for outward minimizing hypersurfaces, which follows from Huisken's and Ilmanen's weak inverse mean curvature flow \cite{HuiskenIlmanen:/2001} or also from \cite{AgostinianiFogagnoloMazzieri:06/2019}. See \cite{McCormick:09/2018,Wei:04/2018} for extensions of this approach to some asymptotically flat manifolds. In the other spaceforms, including de Sitter space, lower bounds for $W_{2}(M)$ were given in \cite{GallegoSolanes:/2005,MakowskiScheuer:11/2016,Natario:08/2015,Scheuer:09/2019,WangXia:07/2014,WeiXiong:/2015}. There are also many results on Minkowski type inequalities with weights, where the mean curvature is integrated against a weight which mostly comes from the ambient geometry; the main candidate is $\vt'$. Such inequalities play a role for the Penrose inequality in general relativity. See \cite{BrendleHungWang:01/2016,ChenLiZhou:04/2019,De-LimaGirao:04/2016,GeWangWu:10/2015,GeWangWuXia:03/2015,KwongMiao:/2014,ScheuerXia:11/2019,Wang:/2015} for various results in this direction.
To the best of my knowledge, there are no Minkowski type inequalities in situations where the ambient space is not asymptotically of constant curvature and hence
 \Cref{thm:mink,,thm:HK,,gen-mink} appear to provide the first such inequalities.
 
In the next section we justify the use of the proposed curvature flow by proving its crucial monotonicity properties. In \Cref{sec:est} we prove a priori estimates for the flow which lead to its convergence in \Cref{sec:proof}. At last, in \Cref{sec:GI}, the proof of the geometric inequalities is completed.

\section{Monotonicity}\label{sec:Moni}

Let $N$ and $M$ be as in \Cref{N} and \Cref{M}.
For the geometric quantities of $N$ and $M$ we use exactly the same notation as in \cite{Scheuer:02/2019}. Hence we do not repeat that part in detail but introduce the most important new objects on the fly.

Additionally we need some Minkowski type formulas, which we deduce here quickly. Denote by $H_{k}$, $1\leq k\leq n$, the normalized $k$-th elementary symmetric polynomial of the principal curvatures $\ka=(\ka_{i})$ of $M$,
\eq{H_{k}=\fr{1}{\binom{n}{k}}\sum_{1\leq i_{1}<\dots< i_{k}\leq n}\ka_{i_{1}}\cdots\ka_{i_{k}}.} 
Furthermore we denote by $H$ the trace of the second fundamental form, i.e.
\eq{H=n H_{1}.}
Let $\Th$ be a primitive of $\vt$. Then we use \cite[equ.~(2.12)]{Scheuer:02/2019} and get
\eq{\label{Theta}\Th_{;ij}=\vt'u_{;i}u_{;j}+\vt u_{;ij}=\vt'u_{;i}u_{;j}+\vt'\vt^{2}\si_{ij}-sh_{ij}=\vt'g_{ij}-sh_{ij},}
where $(g_{ij})$ is the induced metric on $M$ and where we use the outward pointing normal $\nu$ to define the generalized support function
\eq{s=\bar g(\vt\del_{r},\nu)=\bar g(\bar\n\Th,\nu)>0.}
 A semi-colon denotes covariant differentiation with respect the Levi-Civita connection $\n$ of $(g_{ij})$.
The tensor $(h_{ij})$ is the second fundamental form with respect to $-\nu$ and the principal curvatures $(\ka_{i})$ are the eigenvalues of the Weingarten operator
\eq{h^{i}_{j}=g^{ik}h_{ij}.} 
In general we use $(g_{ij})$ to raise and lower indices of tensors. Taking the trace of \eqref{Theta} yields
 \eq{\label{Mink}\int_{M}sH_{1}=\int_{M}\vt'.}
We get a similar relation for $H_{2}$, cf. \cite[Lemma~2.5]{GuanLiWang:/2019} with different notation,
\eq{\label{Mink-2}\int_{M}sH_{2}=\int_{M}\vt'H_{1}-\fr{1}{n(n-1)}\int_{M}\ov{\Rc}(\nu,\n\Th).}
For functions $F$ as in \Cref{F} we use the standard theory of curvature functions $F$ with the conventions as in \cite{Scheuer:02/2019}. In particular, $F$ can be understood to depend on the Weingarten operator $h^{i}_{j}$ or on the two bilinear forms $g_{ij}$ and $h_{ij}$,
\eq{F=F(h^{i}_{j})=F(g_{ij},h_{ij}).}
Then we define
\eq{F^{ij}=\fr{\del F}{\del h_{ij}},\q F^{ij,kl}=\fr{\del^{2}F}{\del h_{ij}\del h_{kl}}.}
See \cite{Andrews:/2007} and \cite[Ch.~2]{Gerhardt:/2006} for more on the theory of curvature functions. In general, latin indices will always range between $1$ and $n$, while we use the Einstein summation convention.

The following lemma is the key to the monotonicity properties required to deduce the geometric inequalities in \Cref{thm:mink,,thm:HK}.				

\begin{lemma}\label{monotone}
Under the assumptions of \Cref{thm:mink}, along \eqref{flow} with $F=\tfr{H_{2}}{H_{1}}$ there hold
\eq{\del_{t}\abs{M_{t}}\geq 0}
with equality for all $t$ precisely if all $M_{t}$ are umbilic,
and
\eq{\label{monotone-5}\del_{t}\br{\int_{M_{t}}H_{1}+\fr 1n\int_{\hat M_{t}}\ov{\Rc}(\del_{r},\del_{r})}\leq 0.}
Under the assumptions of \Cref{thm:HK} there holds
\eq{\del_{t}\abs{\hat M_{t}}\geq 0,}
where equality for all $t$ implies that all $M_{t}$ are umbilic.
\end{lemma}

\pf{

We use the well-known evolution equations, see \cite[Sec~2.3]{Gerhardt:/2006},
\eq{\del_{t}g_{ij}=2\br{\fr{\vt'}{F}-s}h_{ij},}
\eq{\del_{t}\rt{\det g}=\br{\fr{\vt'}{F}-s}H\rt{\det g}}
and 
\eq{\label{ev-h-1}\del_{t}h^{i}_{j}={\br{s-\fr{\vt'}{F}}_{;j}}^{i}+\br{s-\fr{\vt'}{F}}h^{i}_{k}h^{k}_{j}+\br{s-\fr{\vt'}{F}}\ov{\Rm}(x_{;k},\nu,\nu,x_{;j})g^{ki}.}
First we calculate
\eq{\del_{t}\abs{M_{t}}=n\int_{M_{t}}\br{\fr{\vt'H_{1}^{2}}{H_{2}}-sH_{1}}\geq n\int_{M_{t}}\vt'-n\int_{M_{t}}\vt'=0,}
where we have used the Newton-Maclaurin inequality
\eq{H_{2}\leq H_{1}^{2}}
and \eqref{Mink}.
By the equality characterization of the Newton-Maclaurin inequality, in case of equality all $M_{t}$ must be umbilic.
To prove the second claim, we use the decomposition
\eq{\bar\n\Th=\n\Th+s\nu}
and a formula for the Ricci tensor \cite[Prop~2.1]{Brendle:06/2013}
\eq{\ov{\Rc}(\del_{r},\del_{r})=-n\fr{\vt''}{\vt},\q \ov{\Rc}(\del_{r},e_{i})=0,}
where $(e_{i})$ is an orthonormal frame for $(\cS_{0},\si)$.
We use \eqref{Mink-2} to calculate
\eq{\label{monotone-3}\del_{t}\int_{M_{t}}H_{1}&=\fr{1}{n}\int_{M_{t}}H^{2}\br{\fr{\vt'}{F}-s}-\fr{1}{n}\int_{M_{t}}\abs{A}^{2}\br{\fr{\vt'}{F}-s}\\
				&\hp{=}-\fr{1}{n}\int_{M_{t}}\ov{\Rc}(\nu,\nu)\br{\fr{\vt'}{F}-s}\\
				&=(n-1)\int_{M_{t}} H_{2}\br{\fr{\vt'}{F}-s}-\fr 1n\int_{M_{t}}\ov{\Rc}(\nu,\nu)\br{\fr{\vt'}{F}-s}\\
				&=\fr{1}{n}\int_{M_{t}}\ov{\Rc}(\nu,\n\Th)-\fr 1n\int_{M_{t}}\ov{\Rc}(\nu,\nu)\br{\fr{\vt'}{F}-s}\\
				&=\fr 1n\int_{M_{t}}\ov{\Rc}(\nu,\bar\n\Th)-\fr 1n\int_{M_{t}}\ov{\Rc}(\nu,\nu)\fr{\vt'}{F}\\
				&=-\int_{M_{t}}\fr{\vt''}{\vt} s-\fr 1n\int_{M_{t}}\ov{\Rc}(\nu,\nu)\fr{\vt'}{F}.}
Now decompose
\eq{\nu=V+\bar g\br{\nu,\del_{r}}\del_{r}=V+\fr{s}{\vt}\del_{r},}
where $V$ is the projection of $\nu$ onto $\del_{r}^{\perp}$. In the following estimate we first use the representation of $\ov\Rc$ in terms of $\wh\Rc$, \cite[p.~253]{Brendle:06/2013}, and then \eqref{thm:mink-2}:
\eq{\label{monotone-2}\fr{1}{n}\ov{\Rc}(\nu,\nu)&=\fr 1n\ov{\Rc}(V,V)-\fr{\vt''}{\vt}\fr{s^{2}}{\vt^{2}}\\
			&=\fr{1}{n}\br{\wh\Rc(V,V)-\br{\fr{\vt''}{\vt}+(n-1)\fr{\vt'^{2}}{\vt^{2}}}\bar g(V,V)}-\fr{\vt''}{\vt}\fr{s^{2}}{\vt^{2}}\\
			&\geq -\fr{\vt''}{\vt}\bar g(V,V)-\fr{\vt''}{\vt}\fr{s^{2}}{\vt^{2}},}
with equality precisely if \eqref{thm:mink-2} evaluated at $V$ holds with equality or $V=0$.
We have 
\eq{V=\fr{1}{v}(0,\vt^{-2}\si^{ij}u_{;j}),}
where 
\eq{\label{v}v^{2}=1+\vt^{-2}\si^{ij}u_{;i}u_{;j}=\fr{\vt^{2}}{s^{2}},}
and hence
\eq{\bar g(V,V)=\fr{v^{2}-1}{v^{2}}.}
Inserting this into \eqref{monotone-2} gives
\eq{\label{monotone-4}\fr{1}{n}\ov{\Rc}(\nu,\nu)\geq -\fr{\vt''}{\vt}.  }
We use
\eq{\fr{1}{n}\del_{t}\int_{\hat M_{t}}\ov{\Rc}(\del_{r},\del_{r})=-\del_{t}\int_{\hat M_{t}}\fr{\vt''}{\vt}=\int_{M_{t}}\fr{\vt''}{\vt}\br{s-\fr{\vt'}{F}}}
and finally combine this equality with \eqref{monotone-3} and \eqref{monotone-4}, which gives \eqref{monotone-5}.
To estimate the volume, note that the Heintze-Karcher inequality \eqref{thm:HK-2} implies
\eq{\del_{t}\abs{\hat M_{t}}=\int_{M_{t}}\br{\fr{\vt'}{H_{1}}-s}\geq 0.}
}

\section{A priori estimates}\label{sec:est}

In this section we provide all the a priori estimates that are needed to prove \Cref{main:flow}. The existence of a solution to \eqref{flow} on a maximal time interval $[0,T^{*})$ is standard and a proof can be found in \cite[Ch.~2]{Gerhardt:/2006}. There it is also proven that it suffices to get higher regularity estimates for the radial function $u=u(t,\xi)$ which parametrizes the flow hypersurfaces and satisfies the parabolic equation (here we use $\vt'>0$),
\eq{\label{eq:u}\del_{t}u=\br{\fr{\vt'(u)}{F}-s}v^{-1},}
where $v$ is given by \eqref{v}. Once we have estimates for \eqref{eq:u}, we get estimates for the whole flow $x$ as described in \cite[Sec.~2.5]{Gerhardt:/2006}.

Note that in the following estimates the letter ``$c$'' denotes a generic constant which is allowed to depend only on the data of the problem, i.e. on $N$, $x_{0}$ and $F$ and which may change from line to line.
We start with the estimates up to first order.

\begin{lemma}\label{C1-bound}
Under the assumptions of \Cref{main:flow} we have the following estimates for \eqref{flow}.
\enum{
\item \eq{\min_{M_{0}}u(0,\cdot)\leq u(t,\xi)\leq \max_{M_{0}}u(0,\cdot)\q \fa (t,\xi)\in [0,T^{*})\x M_{0},}
\item \eq{v(t,\xi)\leq c\q \fa (t,\xi)\in [0,T^{*})\x M_{0}. }
}

\end{lemma}

\pf{
We use \eqref{Theta} to note that at a maximal point of $u$ there holds
\eq{sh\geq \vt'g}
as bilinear forms. Due to the monotonicity and homogeneity of $F$ we obtain
\eq{F(h^{i}_{j})\geq \fr{\vt'}{s}} 
and thus the function $\max u(t,\cdot)$ is non-increasing in $t$. A similar argument applies to $\min u(t,\cdot)$ and hence the first claim is true.

As $F$ vanishes on the boundary of $\Ga_{+}$, the flows preserves the convexity of $M_{t}$ up to $T^{*}$. This means we have a convex graph in a Riemannian warped product space and the claimed $C^{1}$-estimate is immediate from \cite[Thm.~2.7.10]{Gerhardt:/2006}. 
}

Since \eqref{eq:u} is a fully nonlinear parabolic equation, gradient estimates are not enough to bootstrap up the regularity, as for example in \cite{GuanLi:/2015, GuanLiWang:/2019}. The crucial part is the bound on the curvature. To get this bound, we have to investigate the evolution of $h$, \eqref{ev-h-1}, in greater detail and combine it with other quantities.
Let us define the parabolic operator
\eq{\cP=\del_{t}-\fr{\vt'}{F^{2}}F^{kl}\n_{kl}^{2}-\bar g(\bar\n\Th,\n),}
which may act on functions as well as on time-dependent tensor fields. Here we have to note that we only use time-independent local frames.

We first estimate the curvature function and therefore collect some evolution equations.

\begin{lemma}
Along the flow \eqref{flow} there hold the following evolution equations.
\eq{\label{ev-s}\cP s&=\fr{\vt'}{F^{2}}(F^{ij}h_{ik}h^{k}_{j}-F^{2})s-\fr{1}{F}\bar g(\bar\n\Th,\n\vt')\\
		&\hp{=}+\fr{\vt'}{F^{2}}\bar g(\bar\n\Th,x_{;k})F^{ij}\ov{\Rm}(\nu,x_{;i},x_{;m},x_{;j})g^{mk},}

\eq{\label{ev-thetaprime}\cP \vt'&=\fr{\vt''}{\vt}\fr{2s}{F}\vt'-\vt''\vt-\fr{\vt''}{\vt}\fr{\vt'}{F^{2}}F^{ij}g_{ij}\vt'-\br{\vt'''\vt'-\fr{\vt''\vt'^{2}}{\vt}}\fr{1}{F^{2}}F^{ij}u_{;i}u_{;j}}
and 
\eq{\label{ev-speed}\cP\br{\fr{\vt'}{F}-s}&=\fr{\vt'}{F^{2}}(F^{ij}h_{ik}h^{k}_{j}-F^{2})\br{\fr{\vt'}{F}-s}\\
			&\hp{=}+\fr{\vt'}{F^{2}}F^{ij}\ov{\Rm}(x_{;i},\nu,\nu,x_{;j})\br{\fr{\vt'}{F}-s}+\fr{\vt''}{\vt}\fr{s}{F}\br{\fr{\vt'}{F}-s}.}
\end{lemma}

\pf{
According to \cite[p.~1114]{Scheuer:02/2019} and the Codazzi equation we have
\eq{\label{D2s}s_{;ij}&=\vt'h_{ij}-h_{ik}h^{k}_{j}s+\bar g(\bar\n\Th,x_{;k})h^{k}_{i;j}\\
		&=\vt'h_{ij}-h_{ik}h^{k}_{j}s+\bar g(\bar\n\Th,x_{;k}){h_{ij;}}^{k}-\bar g(\bar\n\Th,x_{;k})\ov{\Rm}(\nu,x_{;i},x_{;m},x_{;j})g^{mk}.}
Also there holds
\eq{\del_{t}s&=\vt'\br{\fr{\vt'}{F}-s}-\bar g\br{\bar\n\Th,\n\br{\fr{\vt'}{F}-s}}\\
		&=\fr{\vt'^{2}}{F}-s\vt'-\fr{1}{F}\bar g(\bar\n\Th,\n\vt')+\fr{\vt'}{F^{2}}\bar g(\bar\n\Th,\n F)+\bar g(\bar\n\Th,\n s).}
Hence the first equation follows.

For the second equation we note
\eq{\del_{t}\vt'=\vt''\del_{t}u=\fr{\vt''}{\vt}\br{\fr{\vt'}{F}-s}s,}
 use \eqref{Theta},
\eq{\label{D2thetaprime}\vt'_{;ij}&=\vt''u_{;ij}+\vt'''u_{;i}u_{;j}\\
			&=\vt''\br{\fr{\vt'}{\vt}g_{ij}-\fr{s}{\vt}h_{ij}-\fr{\vt'}{\vt}u_{;i}u_{;j}}+\vt'''u_{;i}u_{;j}\\
			&=\fr{\vt''\vt'}{\vt}g_{ij}-\fr{\vt''}{\vt}sh_{ij}+\br{\vt'''-\fr{\vt''\vt'}{\vt}}u_{;i}u_{;j}}
and finally note
\eq{-\fr{\vt''}{\vt}s^{2}=-\vt''\vt v^{-2}=-\vt''\vt+\fr{v^{2}-1}{v^{2}}\vt''\vt=-\vt''\vt+\bar{g}(\bar\n\Th,\n \vt').}
			
For the third equation we calculate, using \eqref{ev-h-1},
\eq{\del_{t}\br{\fr{\vt'}{F}-s}&=\fr{\del_{t} \vt'}{F}+\vt'\del_{t}F^{-1}-\del_{t} s\\
				&=\fr{s}{F^{2}}\fr{\vt''\vt'}{\vt}-\fr{\vt''}{\vt}\fr{s^{2}}{F}+\fr{\vt'}{F^{2}}F^{ij}{\br{\fr{\vt'}{F}-s}_{;ij}}+\fr{\vt'}{F^{2}}\br{\fr{\vt'}{F}-s}F^{ij}h_{ik}h^{k}_{j}\\
			&\hp{=}+\fr{\vt'}{F^{2}}\br{\fr{\vt'}{F}-s}F^{ij}\ov{\Rm}(x_{;i},\nu,\nu,x_{;j})\\
			&\hp{=}-\vt'\br{\fr{\vt'}{F}-s}+\bar g\br{\bar\n\Th,\n\br{\fr{\vt'}{F}-s}},	
}
from which the equation follows.
}

\begin{lemma}\label{speed-bound}
Along the flow \eqref{flow} there holds
\eq{\fr{\vt'}{F}-s\leq c.}
\end{lemma}

\pf{
Due to \Cref{C1-bound} there exists $\al>0$ such that
\eq{\vt'\geq 2\al.}
For $\be\in \{0,1\}$ to be determined, the auxiliary function
\eq{w=\log\br{\fr{\vt'}{F}-s}-\log s-\be\log(\vt'-\al)}
is well defined on the open set in spacetime where $\vt'F^{-1}>s$. We prove by maxi\-mum principle that $w$ is bounded from above.  Whenever $w$ is sufficiently large it satisfies the following equation at spacial maximal points:
\eq{\cP w&=\fr{1}{\fr{\vt'}{F}-s}\cP\br{\fr{\vt'}{F}-s}+\fr{\vt'}{F^{2}}F^{kl}\log\br{\fr{\vt'}{F}-s}_{;k}\log\br{\fr{\vt'}{F}-s}_{;l}\\
		&\hp{=}-\fr{1}{s}\cP s-\fr{\vt'}{F^{2}}F^{kl}(\log s)_{;k}(\log s)_{;l}-\fr{\be}{\vt'-\al}\cP \vt'\\
		&\hp{=}-\fr{\be\vt'}{F^{2}}F^{kl}(\log(\vt'-\al))_{;k}(\log(\vt'-\al))_{;l}\\
		&=\fr{1}{\fr{\vt'}{F}-s}\cP\br{\fr{\vt'}{F}-s}-\fr{1}{s}\cP s-\fr{\be}{\vt'-\al}\cP \vt'+\fr{2\be\vt'}{s(\vt'-\al)F^{2}}F^{kl}s_{;k}\vt'_{;l}.}
Due to \eqref{ev-s}, \eqref{ev-thetaprime} and \eqref{ev-speed}, this is
\eq{\label{speed-bound-1}\cP w&=\fr{\vt'}{F^{2}}F^{ij}\ov{\Rm}(x_{;i},\nu,\nu,x_{;j})+\fr{\vt''}{\vt}\fr{s}{F}+\fr{1}{sF}\bar g(\bar\n\Th,\n\vt')\\
		&\hp{=}-\fr{\vt'}{\vt F^{2}}v\bar g(\bar\n\Th,x_{;k})F^{ij}\ov{\Rm}(\nu,x_{;i},x_{;m},x_{;j})g^{mk}-\fr{\vt''}{\vt}\fr{2s}{F}\fr{\be\vt'}{\vt'-\al}\\
		&\hp{=}+\fr{\be}{\vt'-\al}\vt''\vt+\fr{\vt''}{\vt}\fr{\vt'}{F^{2}}F^{ij}g_{ij}\fr{\be\vt'}{\vt'-\al}\\
		&\hp{=}+\fr{\be}{\vt'-\al}\br{\vt'''\vt'-\fr{\vt''\vt'^{2}}{\vt}}\fr{1}{F^{2}}F^{ij}u_{;i}u_{;j}+\fr{2\be\vt'}{s(\vt'-\al)F^{2}}F^{kl}s_{;k}\vt'_{;l}.}
We use
\eq{\bar g(\bar\n\Th,x_{;k})g^{mk}=\vt{u_{;}}^{m},} 
 \cite[equ.~(4.2)]{Scheuer:02/2019} and \cite[p.~1126]{Scheuer:02/2019}:
\eq{\ov{\Rm}(x_{;i},\nu,\nu,x_{;j})&=-\fr{\vt''}{\vt}g_{ij}+\br{\fr{\vt''}{\vt}-\fr{\vt'^{2}}{\vt^{2}}}(\|\n u\|^{2}g_{ij}-u_{;i}u_{;j})+\wt{\Rm}(x_{;i},\nu,\nu,x_{;j})\\
					&=-\fr{\vt''}{\vt}g_{ij}+v\ov{\Rm}(\nu,x_{;i},x_{;m},x_{;j}){u_{;}}^{m},}
where $\wt \Rm$ is the lift of the Riemann tensor of $(\cS_{0},\vt^{2}(r)\si)$ under the projection $\pi\cn  N\ra \cS_{0}$ and its arguments have to be understood as their projections onto $\cS_{0}$.
We get
\eq{&\fr{\vt'}{F^{2}}F^{ij}\ov{\Rm}(x_{;i},\nu,\nu,x_{;j})-\fr{\vt'}{\vt F^{2}}v\bar g(\bar\n\Th,x_{;k})F^{ij}\ov{\Rm}(\nu,x_{;i},x_{;m},x_{;j})g^{mk}\\
	=~&\fr{\vt'}{F^{2}}F^{ij}\ov{\Rm}(x_{;i},\nu,\nu,x_{;j})-\fr{\vt'}{F^{2}}F^{ij}v\ov{\Rm}(\nu,x_{;i},x_{;m},x_{;j}){u_{;}}^{m}\\
	=~&-\fr{\vt''}{\vt}\fr{\vt'}{F^{2}}F^{ij}g_{ij}.}
Returning to \eqref{speed-bound-1} and using $s_{;k}=\vt h^{m}_{k}u_{;m}$ we obtain
\eq{\cP w&\leq c\abs{\vt''}(1+\tfrac{1}{F})-\fr{\vt''}{\vt}\fr{\vt'}{F^{2}}F^{ij}g_{ij}\br{1-\fr{\be\vt'}{\vt'-\al}}\\
	&\hp{=}+\fr{\be}{\vt'-\al}\br{\vt'''\vt'-\fr{\vt''\vt'^{2}}{\vt}}\fr{1}{F^{2}}F^{ij}u_{;i}u_{;j}+\fr{2\be\vt\vt'\vt''}{s(\vt'-\al)F^{2}}F^{kl}h^{m}_{k}u_{;m}u_{;l}.}
Due to \Cref{N}, $\vt''$ can either be globally non-positive or non-negative. In case that $\vt''\geq 0$ we pick $\be=0$ and use the concavity of $F$ which gives $F^{ij}g_{ij}\geq 1$, \cite[Lemma~2.2.19]{Gerhardt:/2006}. At maximal points, where $F$ is very small, we get $\cP w\leq 0$. In case that $\vt''\leq 0$, we pick $\be=1$. Then
\eq{1-\fr{\vt'}{\vt'-\al}=-\fr{\al}{\vt'-\al},}
\eq{\fr{2\vt\vt'\vt''}{s(\vt'-\al)F^{2}}F^{kl}h^{m}_{k}u_{;m}u_{;l}\leq 0,}
since $F^{kl}h^{m}_{k}$ is positive definite due to the convexity of the flow hypersurfaces, and by the assumption on the third derivative of $\vt$ we see that the dominating term has a good sign, which completes this case as well.
}

To complete the a priori estimates we need to prove that the principal curvatures are uniformly bounded from above. For this purpose we deduce the evolution equation of the second fundamental form.

\begin{lemma}\label{lem:ev-h}
Along the flow \eqref{flow}, the Weingarten operator satisfies the following equation.
\eq{\cP h_{j}^{i}&=\fr{\vt'}{F^{2}}F^{kl}h_{lr}h^{r}_{k}h_{j}^{i}-\fr{2\vt'}{F}h^{i}_{k}h_{j}^{k}+\vt'h_{j}^{i}+\fr{1}{F}\fr{\vt''}{\vt}sh_{j}^{i}+s\ov\Rm(x_{;m},\nu,\nu,x_{;j})g^{mi}\\
			&\hp{=}-\fr{1}{F}\fr{\vt''\vt'}{\vt}\de^{i}_{j}+\fr{\vt'}{F^{2}}F^{kl}\ov{\Rm}(\nu,x_{;k},x_{;l},\nu)h_{j}^{i}+2\fr{\vt'}{F}\ov{\Rm}(\nu,x_{;m},\nu,x_{;j})g^{mi}\\
			&\hp{=}-\bar g(\bar\n\Th,x_{;k})\ov{\Rm}(\nu,x_{;l},x_{;m},x_{;j})g^{mk}g^{li}-\fr{1}{F}\br{\vt'''-\fr{\vt''\vt'}{\vt}}{u_{;}}^{i}u_{;j}\\
			&\hp{=}+\fr{{\vt'_{;}}^{i}}{F^{2}}F_{;j}+\fr{\vt'_{;j}}{F^{2}}{F_{;}}^{i}-\fr{2\vt'}{F^{3}}{F_{;}}^{i}F_{;j}+\fr{\vt'}{F^{2}}F^{kl,rs}{h_{kl;}}^{i}h_{rs;j}\\
			&\hp{=}+\fr{\vt'}{F^{2}}F^{kl}\ov{\Rm}(x_{;k},x_{;m},x_{;l},x_{;j})h^{mi}+\fr{\vt'}{F^{2}}F^{kl}\ov{\Rm}(x_{;k},x_{;m},x_{;l},x_{;r})h^m_jg^{ri}\\
			&\hp{=}+2\fr{\vt'}{F^{2}}F^{kl}\ov{\Rm}(x_{;l},x_{;j},x_{;r},x_{;m})h_{k}^mg^{ri}\\
                     &\hp{=} -\fr{\vt'}{F^{2}}F^{kl}\bar\n\ov{\Rm}(\nu,x_{;m},x_{;k},x_{;j},x_{;l})g^{mi}-\fr{\vt'}{F^{2}}F^{kl}\bar\n\ov{\Rm}(\nu,x_{;k},x_{;l},x_{;m},x_{;j})g^{mi}.}
\end{lemma}

\pf{
It is convenient to work with the second fundamental form, which satisfies
\eq{\label{lem:ev-h-1}\del_{t}h_{ij}&=\del_{t}(h^{k}_{i}g_{kj})\\
					&=\br{s-\fr{\vt'}{F}}_{;ij}-\br{s-\fr{\vt'}{F}}h_{ik}h^{k}_{j}+\br{s-\fr{\vt'}{F}}\ov\Rm(x_{;i},\nu,\nu,x_{;j})\\
					&=s_{;ij}-\fr{\vt'_{;ij}}{F}+\fr{\vt'_{;i}}{F^{2}}F_{;j}+\fr{\vt'_{;j}}{F^{2}}F_{;i}+\fr{\vt'}{F^{2}}F_{;ij}-\fr{2\vt'}{F^{3}}F_{;i}F_{;j}\\
					&\hp{=}-\br{s-\fr{\vt'}{F}}h_{ik}h^{k}_{j}+\br{s-\fr{\vt'}{F}}\ov\Rm(x_{;i},\nu,\nu,x_{;j}).}
The main exercise in such calculations is always to turn the term $F_{;ij}$ into a suitable operator on $h_{ij}$. This makes multiple use of the Codazzi-, Gauss- and Weingarten equations and its complexity depends on the restrictions on the ambient space. The main step was already performed in \cite[p.~1111]{Scheuer:02/2019}. From there we conclude that
\eq{F_{;ij}&=F^{kl,rs}h_{kl;i}h_{rs;j}+F^{kl}h_{kl;ij}\\
		&=F^{kl,rs}h_{kl;i}h_{rs;j}+F^{kl}h_{ij;kl}+F^{kl}(h_{la}h_{jk}-h_{lk}h_{ja}+\ov{\Rm}(x_{;l},x_{;j},x_{;k},x_{;a}))h_{i}^a\\  
                     &\hp{=}+F^{kl}(h_{la}h_{ji}-h_{li}h_{ja}+\ov{\Rm}(x_{;l},x_{;j},x_{;i},x_{;a}))h_{k}^a\\
                     &\hp{=}-F^{kl}\bar\n\ov{\Rm}(\nu,x_{;k},x_{;l},x_{;i},x_{;j})-F^{kl}\ov{\Rm}(x_{;m},x_{;k},x_{;l},x_{;i})h^m_j\\
                     &\hp{=}+F^{kl}\ov{\Rm}(\nu,x_{;k},\nu,x_{;i})h_{lj}+F^{kl}\ov{\Rm}(\nu,x_{;k},x_{;l},\nu)h_{ij}\\
                     &\hp{=} -F^{kl}\bar\n\ov{\Rm}(\nu,x_{;i},x_{;k},x_{;j},x_{;l})-F^{kl}\ov{\Rm}(x_{;m},x_{;i},x_{;k},x_{;j})h^m_l\\
                     &\hp{=}+F^{kl}h_{kl}\ov{\Rm}(\nu,x_{;i},\nu,x_{;j})+F^{kl}\ov{\Rm}(\nu,x_{;i},x_{;k},\nu)h_{jl}.}
Now we can use the homogeneity of $F$, that $(F^{i}_{j})$ commutes with $(h^{j}_{k})$ and the symmetries of the curvature tensor to reduce this equation a little bit:

\eq{\label{lem:ev-h-2}F_{;ij}&=F^{kl,rs}h_{kl;i}h_{rs;j}+F^{kl}h_{ij;kl}-Fh_{i}^{a}h_{ja}+F^{kl}h_{la}h^{a}_{k}h_{ij}\\  
                     &\hp{=}+F^{kl}\ov{\Rm}(x_{;k},x_{;m},x_{;l},x_{;j})h_{i}^m+F^{kl}\ov{\Rm}(x_{;k},x_{;m},x_{;l},x_{;i})h^m_j\\
                     &\hp{=}+2F^{kl}\ov{\Rm}(x_{;l},x_{;j},x_{;i},x_{;m})h_{k}^m\\
                     &\hp{=}+F^{kl}\ov{\Rm}(\nu,x_{;k},x_{;l},\nu)h_{ij}+F\ov{\Rm}(\nu,x_{;i},\nu,x_{;j})\\
                     &\hp{=} -F^{kl}\bar\n\ov{\Rm}(\nu,x_{;i},x_{;k},x_{;j},x_{;l})-F^{kl}\bar\n\ov{\Rm}(\nu,x_{;k},x_{;l},x_{;i},x_{;j}).}

We obtain the desired formula by inserting this equation, \eqref{D2s} and \eqref{D2thetaprime} into \eqref{lem:ev-h-1} and reverting to $h^{i}_{j}$ with the help of
\eq{\del_{t}g^{ij}=-2\br{\fr{\vt'}{F}-s}h^{ij}.}		
}

\begin{lemma}\label{curv-bound}
Along the flow \eqref{flow} the principal curvatures are uniformly bounded and range in a compact set of $\Ga_{+}$.
\end{lemma}

\pf{
The proof is similar to the second case in the proof of \cite[Prop.~3.4]{Scheuer:02/2019}. We repeat the main steps for convenience.
As usual, see \cite[Lemma~4.4]{Gerhardt:11/2011} for example, we may define
\eq{w=\log h^{n}_{n}-\log(s-\be)+\al u}
and a bound on $w$ will suffice. Here we work in normal coordinates around a maximum point, 
\eq{g_{ij}=\de_{ij},\q h_{ij}=\ka_{i}\de_{ij}, \q\ka_{1}\leq \dots\leq \ka_{n},}
 $\be$ is small enough and $\al>0$ will be chosen later.
All achieved a priori estimates and \Cref{lem:ev-h} imply
\eq{\cP h^{n}_{n}&\leq \fr{\vt'}{F^{2}}F^{kl}h_{lr}h^{r}_{k}\ka_{n}-\fr{2\vt'}{F}\ka_{n}^{2}+\vt'\ka_{n}+\fr{c}{F}\ka_{n}+c+\fr{c}{F^{2}}F^{ij}g_{ij}(\ka_{n}+1)\\
			&\hp{=}+2\fr{{\vt'_{;n}}}{F^{2}}F_{;n}-\fr{2\vt'}{F^{3}}F_{;n}^{2}+\fr{\vt'}{F^{2}}F^{kl,rs}{h_{kl;n}}h_{rs;n}\\
			&\leq \fr{\vt'}{F^{2}}F^{kl}h_{lr}h^{r}_{k}\ka_{n}-\fr{2\vt'}{F}\ka_{n}^{2}+\vt'\ka_{n}+\fr{c}{F}\ka_{n}+c+\fr{c}{F^{2}}F^{ij}g_{ij}(\ka_{n}+1)\\
			&\hp{=}+\fr{\vt'}{F^{2}}F^{kl,rs}{h_{kl;n}}h_{rs;n}, }
where we used Cauchy-Schwarz.
From \eqref{Theta} we see that
\eq{u_{;ij}=\fr{\vt'}{\vt}g_{ij}-v^{-1}h_{ij}-\fr{\vt'}{\vt}u_{;i}u_{;j}=\fr{\vt'}{\vt}\bar g_{ij}-v^{-1}h_{ij}}
and hence
\eq{\label{ev-u}\cP u&=\br{\tfr{\vt'}{F}-s}v^{-1}-\fr{\vt'}{F^{2}}F^{ij}(\tfr{\vt'}{\vt}\bar g_{ij}-v^{-1}h_{ij})-\vt\|\n u\|^{2}\\
		&=2\fr{\vt'}{F}v^{-1}-\vt-\fr{\vt'}{\vt}\fr{\vt'}{F^{2}}F^{ij}\bar g_{ij}.}
Also using \eqref{ev-s}, we see
\eq{\cP w&\leq -\fr{\be}{s-\be}\fr{\vt'}{F^{2}}F^{kl}h_{lr}h^{r}_{k}-\fr{2\vt'}{F}\ka_{n}+\vt'+\fr{c}{F}(1+\al)+c\ka_{n}^{-1}+\fr{c}{F^{2}}F^{ij}g_{ij}(1+\ka_{n}^{-1})\\
			&\hp{=}+\ka_{n}^{-1}\fr{\vt'}{F^{2}}F^{kl,rs}{h_{kl;n}}h_{rs;n}+\fr{\vt'}{F^{2}}F^{ij}(\log h^{n}_{n})_{;i}(\log h^{n}_{n})_{;j}\\
			&\hp{=}-\fr{\vt'}{F^{2}}F^{ij}(\log(s-\be))_{;i}(\log(s-\be))_{;j}-\al\fr{\vt'}{\vt}\fr{\vt'}{F^{2}}F^{ij}\bar g_{ij}. }
It is necessary to pick $\al$ large. In order to deal with the resulting derivative coming from the replacement of $(\log h^{n}_{n})_{;i}$, we use a trick that was already used in \cite{Enz:10/2008}. The concavity of $F$ and $\ka_{1}>0$ implies for all symmetric matrices $(\eta_{kl})$:
\eq{F^{nn}\leq \dots \leq F^{11}\q\mbox{and}\q F^{kl,rs}\eta_{kl}\eta_{rs}\leq \fr{2}{\ka_{n}}\sum_{k=1}^{n}(F^{nn}-F^{kk})\eta_{nk}^{2}. }
We apply this to $\eta_{kl}=h_{kl;n}$ and estimate
\eq{\ka_{n}^{-1}F^{kl,rs}h_{kl;n}h_{rs;n}&\leq \fr{2}{\ka_{n}^{2}}\sum_{k=1}^{n}(F^{nn}-F^{kk})(h_{kn;n})^{2}\\
					&\leq \fr{2}{\ka_{n}^{2}}\sum_{k=1}^{n}(F^{nn}-F^{kk})(h_{nn;k})^{2}+\fr{c}{\ka_{n}^{2}}\sum_{k=1}^{n}(F^{kk}-F^{nn}),}
due to the Codazzi equation.
Thus, at a maximal point of $w$,
\eq{&F^{ij}(\log h^{n}_{n})_{;i}(\log h^{n}_{n})_{;j}+\ka_{n}^{-1}F^{kl,rs}h_{kl;n}h_{rs;n}\\
	\leq~& \fr{1}{\ka_{n}^{2}}F^{nn}\sum_{k=1}^{n}(h_{nn;k})^{2}+\fr{c}{\ka_{n}^{2}}\sum_{k=1}^{n}(F^{kk}-F^{nn})\\
	=~& F^{nn}\sum_{k=1}^{n}(\log(s-\be)_{;k}-\al u_{;k})^{2}+\fr{c}{\ka_{n}^{2}}\sum_{k=1}^{n}(F^{kk}-F^{nn})\\
	\leq ~&c\ka_{n}^{-2}F^{ij}g_{ij}+F^{nn}\|\n(\log(s-\be))\|^{2}+\al^{2}F^{nn}\|\n u\|^{2}\\
	\hp{=}&-2\al F^{nn}g(\n \log(s-\be),\n u).}
The estimate on $\cP w$ becomes at a spacetime maximum
\eq{0\leq\cP w&\leq -\fr{\be}{s-\be}\fr{\vt'}{F^{2}}F^{kl}h_{lr}h^{r}_{k}-\fr{2\vt'}{F}\ka_{n}+\vt'+\fr{c}{F}(1+\al)+c\ka_{n}^{-1}\\
			&\hp{=}+\fr{c}{F^{2}}F^{ij}g_{ij}(1+\ka_{n}^{-1})+c\ka_{n}^{-2}\fr{\vt'}{F^{2}}F^{ij}g_{ij}+\fr{\vt'}{F^{2}}F^{nn}\|\n(\log(s-\be))\|^{2}\\
			&\hp{=}+\al^{2}\fr{\vt'}{F^{2}}F^{nn}\|\n u\|^{2}-2\al\fr{\vt'}{F^{2}} F^{nn}g(\n \log(s-\be),\n u)\\
			&\hp{=}-\fr{\vt'}{F^{2}}F^{ij}(\log(s-\be))_{;i}(\log(s-\be))_{;j}-\al\fr{\vt'}{\vt}\fr{\vt'}{F^{2}}F^{ij}\bar g_{ij}\\
			&\leq \fr{\vt'}{F^{2}}F^{nn}\br{-\fr{\be}{s-\be}\ka_{n}^{2}+\al^{2}\|\n u\|^{2}+\al c\ka_{n}}-\fr{\vt'}{F}\ka_{n}+\fr{c}{F}(1+\al)\\
			&\hp{=}+\fr{\vt'}{F^{2}}F^{ij}g_{ij}\br{c\ka_{n}^{-2}-\ep_{0}\al\fr{\vt'}{\vt}+c+c\ka_{n}^{-1}}, }
where we also used 
\eq{F\leq \ka_{n}\q\mbox{and}\q F^{ij}\bar g_{ij}\geq \ep_{0}F^{ij}g_{ij}}
for some constant $\ep_{0}$. Fixing a sufficiently large $\al$, we see that $\ka_{n}$ can not be too large without reaching a contradiction. Hence $\ka_{n}$ is bounded. Due to the lower bound on $F$, the principal curvatures range in a compact subset of its domain. The proof is complete.
}

\begin{cor}
The flow \eqref{flow} exists for all times with uniform $C^{\8}$-estimates.
\end{cor}

\pf{
This is standard. The radial function $u$ satisfies a fully nonlinear PDE
\eq{\del_{t}{u}=\br{\fr{\vt'}{F}-s}v^{-1}=G(\cdot,u,\n u,\n^{2}u),}
which is uniformly parabolic due to \cref{speed-bound,curv-bound}. As we assumed that $F$ is a concave curvature function, $G$ is concave in $\n^{2}u$. Furthermore we have $C^{2}$-bounds due to all apriori estimates. We may apply the regularity results by Krylov-Safonov \cite{Krylov:/1987} and linear Schauder theory to obtain uniform bounds on all derivatives of $u$. A standard continuation argument \cite[Sec.~2.5]{Gerhardt:/2006} proves the long-time existence.
}

\section{Completion of the convergence proof}\label{sec:proof}

\begin{proof}[Proof of \Cref{main:flow}]
The proof is similar to the one in \cite[Sec.~5]{LambertScheuer:10/2019}. We sketch the idea. From \eqref{ev-u} we see that $\Th$ satisfies
\eq{\cP\Th=2\fr{\vt'}{F}s-\vt^{2}-\fr{\vt'^{2}}{F^{2}}F^{ij}g_{ij}\leq 0.}
The strong maximum principle implies that 
\eq{\osc\Th(t)=\max_{M_{t}}\Th-\min_{M_{t}}\Th}
is strictly decreasing, unless it is constant, in which case $M_{t}$ is already a  radial slice. Suppose that $\osc\Th$ does not converge to zero, but to a positive value $\al>0$.
From the smooth a priori estimates we can define a smooth limit flow (possibly after choosing a subsequence)
\eq{x_{\8}(t,\xi)=\lim_{i\ra \8}x(t+i,\xi).}
This limit flow has constant oscillation $\al>0$, a contradiction to the strong maximum principle, as $x_{\8}$ satisfies \eqref{flow} as well. 
Hence every subsequential limit of $(M_{t})$ must be a radial slice and due to the barrier estimates it can only be a unique slice.
\end{proof}

\section{Geometric inequalities}\label{sec:GI}
\begin{proof}[Proof of \Cref{thm:mink,,thm:HK}]
With this convergence result at hand, the results in \Cref{thm:mink,,thm:HK} follow immediately from the monotonicity properties in \Cref{monotone} and their equality characterizations.
\end{proof}

\subsection*{Non-positive radial curvature}

We use an idea of Simon Brendle, who used the inverse mean curvature flow \cite{Gerhardt:11/2011} to prove \eqref{gen-mink-1} in the hyperbolic space, see \cite{GuanLi:/2019}. We adapt this proof to the ambient spaces given in \Cref{gen-mink} and use the result on inverse mean curvature flow from \cite{Scheuer:02/2019}.

\begin{proof}[Proof of \Cref{gen-mink}]
Denote by $S_{r}$ the radial $r$-slice in $N$ and write
\eq{W_{2}(M)=\int_{M}H_{1}-\abs{\hat M}.}
Then along the inverse mean curvature flow in $N$,
\eq{\del_{t}x=\fr{1}{H}\nu,}
we have the following variational formulae:
\eq{\del_{t}\int_{M_{t}}H=\int_{M_{t}}H-\int_{M_{t}}\fr{1}{H}\br{\|A\|^{2}+\ov{\Rc}(\nu,\nu)}=2\int_{M_{t}}\fr{\si_{2}}{H}-\int_{M_{t}}\fr{1}{H}\ov{\Rc}(\nu,\nu).}
From \eqref{monotone-4} we obtain
\eq{\del_{t}\int_{M_{t}}H\leq (n-1)\int_{M_{t}}\fr{H_{2}}{H_{1}}+\int_{M_{t}}\fr{n}{H},}
where in case that \eqref{gen-mink-2} holds strictly, this inequality is strict unless we have a flow of slices.

Also we obtain
\eq{\del_{t}\abs{\hat M_{t}}=\int_{M_{t}}\fr{1}{H}.}
Let $\phi$ be defined by the relation
\eq{W_{2}(S_{r})=\phi(\abs{S_{r}}),}
which is well-defined due to the strict monotonicity of $W_{2}(S_{r})$ and $\abs{S_{r}}$ with respect to $r$. Hence along inverse mean curvature flow we have
\eq{0=\del_{t}(W_{2}(S_{r(t)})-\phi(\abs{S_{r(t)}}))&=\fr{n-1}{n}\int_{S_{r(t)}}\fr{H_{2}}{H_{1}}-\phi'(\abs{S_{r(t)}})\abs{S_{r(t)}}\\
					&=\fr{n-1}{n}\int_{S_{r(t)}}H_{1}-\phi'(\abs{S_{r(t)}})\abs{S_{r(t)}}\\
					&=\fr{n-1}{n}\br{W_{2}(S_{r(t)})+\abs{\hat S_{r(t)}}}-\phi'(\abs{S_{r(t)}})\abs{S_{r(t)}}} 
and hence
\eq{\phi'(\abs{S_{r(t)}})\abs{S_{r(t)}}=\fr{n-1}{n}\br{W_{2}(S_{r(t)})+\abs{\hat S_{r(t)}}}.}

Now, for a general mean-convex flow $(M_{t})$ pick a flow of spheres $(S_{r(t)})$ such that
\eq{\abs{M_{t}}=\abs{S_{r(t)}}}
 and calculate
\eq{\del_{t}\br{W_{2}(t)-\phi(\abs{M_{t}})}&\leq \fr{n-1}{n}\int_{M_{t}}\fr{H_{2}}{H_{1}}-\phi'(\abs{M_{t}})\abs{M_{t}}\\
					&\leq \fr{n-1}{n}\int_{M_{t}}H_{1}-\phi'(\abs{M_{t}})\abs{M_{t}}\\
					&= \fr{n-1}{n}\br{W_{2}(t)+\abs{\hat M_{t}}}-\fr{n-1}{n}W_{2}(S_{r}(t))-\fr{n-1}{n}\abs{\hat S_{r(t)}}\\
					&= \fr{n-1}{n}(W_{2}(t)-\phi(\abs{M_{t}}))+\fr{n-1}{n}\br{\abs{\hat M_{t}}-\abs{\hat S_{r(t)}}}\\
					&\leq \fr{n-1}{n}(W_{2}(t)-\phi(\abs{M_{t}})),}
due to the isoperimetric inequality \cite[Thm.~1.2]{GuanLiWang:/2019}.
Hence
\eq{W_{2}(t)-\phi(\abs{M_{t}})\leq e^{\fr{n-1}{n}t}(W_{2}(0)-\phi(\abs{M_{0}})),}
where equality implies total umbilicity.
Now we use that for large $t$ the inverse mean curvature flow is almost umbilic in the sense that
\eq{\left|h^{i}_{j}-\fr{\vt'}{\vt}\de^{i}_{j}\right|\leq \fr{ct}{\vt'\vt},}
see \cite[equ.~(1.3)]{Scheuer:02/2019}.
Since 
\eq{\fr{\sinh(r)}{\cosh(r)}\leq \fr{\vt'}{\vt}\leq \fr{\cosh(r)}{\sinh(r)},}
$\vt'/\vt$ is uniformly bounded above and below by positive constants. As geodesic spheres satisfy the ODE
\eq{\fr{d}{dt}\vt(r)=\fr 1n\vt(r)}
and any solution of inverse mean curvature flow respects initial spherical barriers, we can estimate
\eq{\vt\geq ce^{\fr{t}{n}} }
and similarly for $\vt'$. Hence
\eq{H=n\fr{\vt'}{\vt}+O(te^{-\fr{2t}{n}})}
and we estimate with the help of the isoperimetric inequality,
\eq{W_{2}(t)&=\fr{1}{n}\int_{M_{t}}H-\abs{\hat M_{t}}\\
		&=\int_{M_{t}}\fr{\vt'}{\vt}+O(te^{-\fr{2t}{n}})\abs{M_{t}}-\abs{\hat M_{t}}\\
		&=\int_{M_{t}}\br{\fr{\vt'}{\vt}-\fr{\vt'(r(t))}{\vt(r(t))}}+\fr{\vt'(r(t))}{\vt(r(t))}\abs{S_{r(t)}}+O(te^{-\fr{2t}{n}})\abs{M_{t}}-\abs{\hat M_{t}}\\
		&\geq W_{2}(S_{r(t)})+\int_{M_{t}}\br{\fr{\vt'}{\vt}-\fr{\vt'(r(t))}{\vt(r(t))}}+O(te^{-\fr{2t}{n}})\abs{M_{t}}.}
Now we have to estimate the integral term. There holds
\eq{\left|\fr{\vt'(r(t))}{\vt(r(t))}-1\right|=\fr{1}{\vt(r(t))}\fr{\abs{\vt'(r(t))^{2}-\vt(r(t))^{2}}}{\abs{\vt'(r(t))+\vt(r(t))}}=\fr{(\al^{2}-\be^{2})\vt(r(t))^{-1}}{\abs{\vt'(r(t))+\vt(r(t))}}=O(e^{-\fr{2t}{n}})}
and similarly for $\vt'/\vt$. Hence
\eq{W_{2}(t)\geq W_{2}(S_{r(t)})+O(te^{-\fr{2t}{n}})\abs{M_{t}}.}
Finally,
\eq{W_{2}(t)-\phi(\abs{M_{t}})\geq W_{2}(S_{r(t)})-\phi(\abs{S_{r(t)}})+O(te^{-\fr{2t}{n}})\abs{M_{t}}=O(te^{-\fr{2t}{n}})\abs{M_{t}}}
and hence
\eq{W_{2}(0)-\phi(\abs{M_{0}})\geq e^{-\fr{n-1}{n}t}O(te^{-\fr{2t}{n}})e^{t}\abs{M_{0}}=O(te^{-\fr{t}{n}})\abs{M_{0}}\ra 0,\q t\ra \8. }
This completes the proof.
\end{proof}

\section*{Acknowledgments}
This work was made possible through a research scholarship the author received from the DFG and which was carried out at Columbia University in New York. JS would like to thank the DFG, Columbia University and especially Prof.~Simon Brendle for their support.

\bibliographystyle{/Users/J_Mac/Documents/Uni/TexTemplates/shamsplain}
\bibliography{/Users/J_Mac/Documents/Uni/TexTemplates/Bibliography.bib}

\end{document}